\pdfoutput=1
\RequirePackage{ifpdf}
\ifpdf 
\documentclass[pdftex]{sigma}
\else
\documentclass{sigma}
\fi

\numberwithin{equation}{section}

\def\eps{\varepsilon}
\def\dsp{\displaystyle}
\def\bigO{{\cal O}}
\def\calC{{{\cal C}}}
\def\calL{{{\cal L}}}

\begin{document}
\allowdisplaybreaks

\newcommand{\arXivNumber}{1804.06749}

\renewcommand{\thefootnote}{}

\renewcommand{\PaperNumber}{073}

\FirstPageHeading

\ShortArticleName{Expansions of Jacobi Polynomials for Large Values of $\beta$ and of Their Zeros}

\ArticleName{Asymptotic Expansions of Jacobi Polynomials\\ for Large Values of $\boldsymbol{\beta}$ and of Their Zeros\footnote{This paper is a~contribution to the Special Issue on Orthogonal Polynomials, Special Functions and Applications (OPSFA14). The full collection is available at \href{https://www.emis.de/journals/SIGMA/OPSFA2017.html}{https://www.emis.de/journals/SIGMA/OPSFA2017.html}}}

\Author{Amparo GIL~$^\dag$, Javier SEGURA~$^\dag$ and Nico M.~TEMME~$^\S$}

\AuthorNameForHeading{A.~Gil, J.~Segura and N.M.~Temme}

\Address{$^\dag$~Departamento de Matem\'atica Aplicada y CC, de la Computaci\'on, ETSI Caminos, \\
\hphantom{$^\dag$}~Universidad de Cantabria, 39005-Santander, Spain}
\EmailD{\href{mailto:amparo.gil@unican.es}{amparo.gil@unican.es}}
\URLaddressD{\url{http://personales.unican.es/gila/}}

\Address{$^\ddag$~Departamento de Matem\'aticas, Estadistica y Computaci\'on,\\
\hphantom{$^\ddag$}~Universidad de Cantabria, 39005 Santander, Spain}
\EmailD{\href{mailto:javier.segura@unican.es}{javier.segura@unican.es}}

\Address{$^\S$~IAA, 1825 BD 25, Alkmaar, The Netherlands\\
\hphantom{$^\S$}~Former address: Centrum Wiskunde \& Informatica (CWI),\\
\hphantom{$^\S$}~Science Park 123, 1098 XG Amsterdam, The Netherlands}
\EmailD{\href{mailto:nico.temme@cwi.nl}{nico.temme@cwi.nl}}

\ArticleDates{Received April 19, 2018, in final form July 12, 2018; Published online July 17, 2018}

\Abstract{Asymptotic approximations of Jacobi polynomials are given for large values of the $\beta$-parameter and of their zeros. The expansions are given in terms of Laguerre polynomials and of their zeros. The levels of accuracy of the approximations are verified by numerical examples.}


\Classification{33C45; 41A60; 65D20}

\renewcommand{\thefootnote}{\arabic{footnote}}
\setcounter{footnote}{0}

\section{Introduction}

This paper fits in a series of our papers on asymptotic aspects of Gauss quadrature for the classical polynomials. In \cite{Gil:2018:GHL} we have considered large degree asymptotics of Hermite and Laguerre polynomials, including asymptotic expansions of their nodes and weights for Gauss quadrature. In a recently submitted paper \cite{Gil:2018:NIC} we study similar aspects for large degree asymptotics of the Jacobi polynomials and in the present paper we do the same for large values of~$\beta$.

We describe two rather similar methods, the second one gives more concise coefficients and better approximations. Our methods are different from the large $\beta$ approach used in \cite{Dimitrov:2016:ABJ}, where the starting point is a power series representation of the Jacobi polynomial. We start with an integral representation.

Similar expansions can be obtained for large values of $\alpha$ with other parameters fixed. This follows from the symmetry relation
\begin{gather*}
P_n^{(\alpha,\beta)}(x)=(-1)^n P_n^{(\beta,\alpha)}(-x).
\end{gather*}

\section{A first expansion}\label{sec:Jacplargebeta}
For details and properties of the orthogonal polynomials, we refer to \cite{Koornwinder:2010:OPS}. The following limit is well known
\begin{gather*}
\lim_{\beta\to\infty} P_n^{(\alpha,\beta)}\left(1-\frac{2x}{\beta}\right)=L_n^{(\alpha)}(x).
\end{gather*}
The expansions in this paper give asymptotic details of this limit.

For deriving the results we use the Cauchy integral that follows from the Rodrigues formula for the Jacobi polynomial
\begin{gather}\label{eq:Jacbeta02}
P_n^{(\alpha,\beta)}(x) =\frac{(-1)^n}{2^n w(x)}\frac{1}{2\pi i}\int_\calC \frac{w(z)\big(1-z^2\big)^n}{(z-x)^{n+1}} {\rm d}z, \qquad x\in(-1,1),
\end{gather}
where $w(z)=(1-z)^{\alpha}(1+z)^{\beta}$ and $\calC$ is a contour around $x$ with $\pm1$ outside the contour. Similarly, for the Laguerre polynomial
\begin{gather*}
L_n^{(\alpha)}(x) =\frac{1}{2\pi i}\int_\calL (s+1)^{n+\alpha}e^{-xs} \frac{{\rm d}s}{s^{n+1}},
\end{gather*}
where $\calL$ is a contour around $0$ with $-1$ outside the contour.

By changing the variable in \eqref{eq:Jacbeta02}, writing $z=x-(1-x)s$, it follows that
\begin{gather}\label{eq:Jacbeta04}
P_n^{(\alpha,\beta)}(x)= \frac{(1+x)^n}{2^n}\frac{1}{2\pi i}\int_\calC (s+1)^{n+\alpha}\left(1-\frac{1-x}{1+x}s\right)^{n+\beta}\frac{ {\rm d}s}{s^{n+1}},
\end{gather}
where $\calC$ is a contour around $0$ with $-1$ and $(1+x)/(1-x)$ outside the contour.

Next, replacing $x$ by $1-2x/\beta$ we find for $0<x<\beta$
\begin{gather}\label{eq:Jacbeta05}
P_n^{(\alpha,\beta)}\left(1-\frac{2x}{\beta}\right)=\frac{(1-x/\beta)^n}{2\pi i}\int_\calC (s+1)^{n+\alpha}\left(1-\frac{x}{\beta-x}s\right)^{n+\beta}\frac{ {\rm d}s}{s^{n+1}},
\end{gather}
where the contour is around $0$ with $-1$ and $(\beta-x)/x$ outside the contour.

We expand
\begin{gather}\label{eq:Jacbeta06}
(1-x/\beta)^n\left(1-\frac{x}{\beta-x}s\right)^{n+\beta}=e^{-xs}\sum_{k=0}^\infty\frac{c_k(n,x;s)}{\beta^k},
\end{gather}
and we obtain
\begin{gather}\label{eq:Jacbeta07}
P_n^{(\alpha,\beta)}\left(1-\frac{2x}{\beta}\right)=\sum_{k=0}^n\frac{\Phi_k(n,\alpha,x)}{\beta^k},
\end{gather}
where
\begin{gather}\label{eq:Jacbeta08}
\Phi_k(n,\alpha,x)= \frac{1}{2\pi i}\int_\calC (s+1)^{n+\alpha}e^{-xs}c_k(n,\alpha,x;s)\frac{ {\rm d}s}{s^{n+1}}.
\end{gather}

The same expansion has been derived in \cite{Dimitrov:2016:ABJ}, starting from a polynomial ${}_2F_1$-representation of the Jacobi polynomial, and explicit forms of $\Phi_k(n,\alpha,x)$ are given in terms of symmetric functions. The first $\Phi_k(n,\alpha,x)$ are shown in terms of derivatives of $L_{n}^{(\alpha)}(x)$.

\begin{remark}\label{rem:rem01}
The expansion in \eqref{eq:Jacbeta06} is convergent for $\beta> \vert x\vert\left(1+\vert s\vert\right)$, and initially we obtain an infinite expansion in \eqref{eq:Jacbeta07} of the Jacobi polynomial. Because this expansion defines an analytic function of $\beta$ for large enough values of $\beta$ (depending on $x$ and $s$), the relation in \eqref{eq:Jacbeta07} (with the infinite series) defines an analytical identity. But when we expand the left-hand side into an $n$th degree power series by using one of the available power series of a ${}_2F_1$-representation, we can rearrange that expansion into a finite $n$th degree expansion with negative powers of $\beta$. Because of the analytical relationship the series in~\eqref{eq:Jacbeta07} is finite, as shown.

Another point is that the Rodrigues formula and the Cauchy integral in~\eqref{eq:Jacbeta02} are not only valid for $x\in(-1,1)$, as mentioned because of convention. But the relation in~\eqref{eq:Jacbeta07} shows two polynomials in $x$ on both sides. Again, by using the analytical relationship we can take for~$x$ any complex value. But the large $\beta$ asymptotic property holds uniformly for bounded values of~$n$,~$\alpha$ and~$x$.
\end{remark}

The first coefficients $c_k(n,x;s)$ are
\begin{gather}
c_0(n,x;s)=1,\nonumber\\
c_1(n,x;s) = -\tfrac12x\bigl(2xs+xs^2+2ns+2n\bigr),\nonumber\\
c_2(n,x;s) = \tfrac{1}{24}x^2\bigl({-}24xs^2-24xs-8xs^3-24ns-12ns^2 + 12x^2s^2+12x^2s^3+36xs^2n\nonumber\\
\hphantom{c_2(n,x;s) =}{} +3x^2s^4+12xs^3n +12n^2s^2+24nxs+24n^2s-12n+12n^2\bigr).\label{eq:Jacbeta09}
\end{gather}
These are polynomials in $s$ of degree $2k$. The $\Phi_k(n,\alpha,x)$ are combinations of Laguerre polynomials, because a power $s^j$ in $c_k(n,x;s)$ gives
\begin{gather}\label{eq:Jacbeta10}
L_{n-j}^{(\alpha+j)}(x)=(-1)^j\frac{ {\rm d}^j}{ {\rm d}x^j}L_{n}^{(\alpha)}(x).
\end{gather}
That is, when we write
\begin{gather}\label{eq:Jacbeta11}
c_k(n,x;s)=\sum_{j=0}^{2k}c_{jk}s^j,
\end{gather}
we obtain
\begin{gather}\label{eq:Jacbeta12}
\Phi_k(n,\alpha,x)=\sum_{j=0}^{\min(n,2k)}c_{jk}L_{n-j}^{(\alpha+j)}(x)=\sum_{j=0}^{\min(n,2k)}c_{jk} (-1)^j\frac{ {\rm d}^j}{ {\rm d}x^j}L_{n}^{(\alpha)}(x).
\end{gather}

For expansions of the zeros of the Jacobi polynomial when $\beta$ is large, it is convenient to have a representation in the form
\begin{gather}\label{eq:Jacbeta13}
P_n^{(\alpha,\beta)}\left(1-\frac{2x}{\beta}\right)=L_{n}^{(\alpha)}(x)U(n,\alpha,\beta,x)+L_{n-1}^{(\alpha)}(x)V(n,\alpha,\beta,x),
\end{gather}
with expansions
\begin{gather}\label{eq:Jacbeta14}
U(n,\alpha,\beta,x)=\sum_{k=0}^n\frac{u_k(n,\alpha,x)}{\beta^k},\qquad V(n,\alpha,\beta,x)=\sum_{k=0}^n\frac{v_k(n,\alpha,x)}{\beta^k}.
\end{gather}
These coefficients can be obtained by writing the Laguerre polynomials $L_{n-j}^{(\alpha+j)}(x)$ in~\eqref{eq:Jacbeta10} as (see also \cite[Lemma~1]{Dimitrov:2016:ABJ})
\begin{gather*}
L_{n-j}^{(\alpha+j)}(x)=P_j(n,\alpha,x)L_{n}^{(\alpha)}(x)+Q_j(n,\alpha,x)L_{n-1}^{(\alpha)}(x),
\end{gather*}
where $j=0,1,2,\dots,n$ and the $P_j$ and $Q_j$ follow from
\begin{gather*}
P_0(n,\alpha,x)=1,\qquad Q_0(n,\alpha,x)=0,\qquad
xP_1(n,\alpha,x)=-n,\qquad xQ_1(n,\alpha,x)=n+\alpha,\\
xP_{j+1}(n,\alpha,x)=(\alpha+j-x)P_j(n,\alpha,x)+(j-n-1)P_{j-1}(n,\alpha,x),\\
xQ_{j+1}(n,\alpha,x)=(\alpha+j-x)Q_j(n,\alpha,x)+(j-n-1)Q_{j-1}(n,\alpha,x).
\end{gather*}
These relations follow from the standard recurrence relations of the Laguerre polynomials. An extra useful relation is (for $j= 0,1,\dots,n-1$),
\begin{gather*}
xL_{n-j-1}^{(\alpha+j+1)}(x)+(x-j-\alpha)L_{n-j}^{(\alpha+j)}(x)+(n-j+1)L_{n-j+1}^{(\alpha+j-1)}(x)=0.
\end{gather*}

The coefficients $u_k(n,\alpha,x)$ and $v_k(n,\alpha,x)$ are given by
\begin{gather}\label{eq:Jacbeta18}
u_k(n,\alpha,x)=\sum_{j=0}^{\min(n,2k)}c_{jk}P_j(n,\alpha,x),\qquad v_k(n,\alpha,x) = \sum_{j=0}^{\min(n,2k)}c_{jk}Q_j(n,\alpha,x).
\end{gather}
The first coefficients are
\begin{gather*}
u_0(n,\alpha,x)=1,\qquad v_0(n,\alpha,x)=0,\\
u_1(n,\alpha,x) = \tfrac12n(2n+\alpha+1),\qquad v_1(n,\alpha,x) = -\tfrac12(n+\alpha)(\alpha+2n+x+1).
\end{gather*}

We give a numerical example. We take $n=10$, $\alpha=\frac13$, $x=1$, and several values of $\beta$. We use the expansions in \eqref{eq:Jacbeta14} summing up till $k=k_{\max}$. In Table~\ref{tab:tab01} we give the relative errors. We compared the results with computed results by Maple with Digits = 16.

\begin{remark}\label{rem:rem02}
The results become better when we take smaller values of $x$. First observe that the coefficients $c_k(n,\alpha,x;s)$ used in \eqref{eq:Jacbeta06} and~\eqref{eq:Jacbeta08} have a front factor $x^k$, a few first examples are shown in \eqref{eq:Jacbeta09}. Hence, for $x=0$ only the term $k=0$ survives and $x=0$ gives on both sides of~\eqref{eq:Jacbeta07} the correct value
\begin{gather*}
P_n^{(\alpha,\beta)}(1)=L_{n}^{(\alpha)}(0)=\binom{n+\alpha}{n}.
\end{gather*}
This property of the coefficients $c_k(n,\alpha,x;s)$ is lost in the results \eqref{eq:Jacbeta13}, which is obtained after rearranging \eqref{eq:Jacbeta07}. However, using $x=0$ in \eqref{eq:Jacbeta13} with $k=0$, that is, with $U=1$ and $V=0$, gives again the correct result, but $U$ and $V$ assume different values when using the sums in~\eqref{eq:Jacbeta14}.

As an example that small values of $x$ give better results, we use $x=1/100$. With $\beta=50$, $n=10$, $\alpha=1/3$ as in Table~\ref{tab:tab01}, we obtain the results for various values of $k_{\max}$, which are much better than those for $x=1$:
\begin{gather*}
0.86\times 10^{-4}, \qquad 0.19\times 10^{-6}, \qquad 0.22\times 10^{-9}, \qquad 0.15e\times 10^{-12}, \qquad 0.90\times 10^{-15}.
\end{gather*}
When we want to use the expansion of $P_n^{(\alpha,\beta)}(z)$ for a general $z\in(-1,1)$, the correspon\-ding~$x$ in $P_n^{(\alpha,\beta)}(1-2x/\beta)$ is $x=\frac12\beta(1-z)$. The representation in \eqref{eq:Jacbeta13} and expansions in \eqref{eq:Jacbeta14} remain valid, but the expansions lose their asymptotic property, because $x=\bigO(\beta)$.
\end{remark}

\begin{table}[h]
\caption{Relative errors in the computation of $P_n^{(\alpha,\beta)}(1-2x/\beta)$ with $n=10$, $\alpha=\frac13$, $x=1$, and several values of $\beta$ by using the expansions \eqref{eq:Jacbeta14} summing up till $k=k_{\max}$.}\label{tab:tab01}
$$
\begin{array}{cccccc}
k_{\max} \to & 1 & 2&3&4 &5 \\
\beta &&&&&\\
\hline
 \hphantom{00}50 \quad & 0.19\times 10^{-0} & 0.13 \times 10^{-1} & 0.12 \times 10^{-1} & 0.16 \times 10^{-2} & 0.99 \times 10^{-4}\tsep{2pt}\\
 \hphantom{0}100 \quad& 0.67 \times 10^{-1} & 0.28 \times 10^{-2} & 0.94 \times 10^{-3} & 0.60 \times 10^{-4} & 0.19 \times 10^{-5}\\
 \hphantom{0}500 \quad& 0.38 \times 10^{-2} & 0.39 \times 10^{-4} & 0.22 \times 10^{-5} & 0.27 \times 10^{-7} & 0.17 \times 10^{-9}\\
1000 \quad & 0.10 \times 10^{-2} & 0.53 \times 10^{-5} & 0.14 \times 10^{-6} & 0.90 \times 10^{-9} & 0.28 \times 10^{-11}\\
\hline
\end{array}
$$
\end{table}

\subsection{Expansions of the zeros}\label{sec:Jacnlargebetazer}
We use the method described in our previous papers \cite{Gil:2018:GHL,Gil:2018:NIC}. We denote the zeros of
\begin{alignat*}{3}
&P_n^{(\alpha,\beta)}(z) \qquad && {\rm by}\quad z_k, \quad z_1<z_2<\cdots < z_n,&\\
&P_n^{(\alpha,\beta)} (1-2x/\beta )\qquad && {\rm by} \quad x_k, \quad x_1>x_2>\cdots >x_n,& \\
&L_n^{(\alpha)}(x)\qquad & &{\rm by} \quad \ell_k, \quad \ell_1<\ell_2<\cdots < \ell_n.&
\end{alignat*}
Clearly, for large $\beta$, $x_1$ can be approximated by $\ell_n$, $x_2$ by $\ell_{n-1}$, and in general $x_k$ by $\ell_{n-k+1}$, $k=1,2,\dots,n$.

We assume for $x_k$ an expansion of the form
\begin{gather}\label{eq:Jacbeta23}
x_k=\ell_{n-k+1}+\eps,\qquad \eps\sim\sum_{k=1}^\infty \frac{\eps_k}{\beta^k},\qquad \beta\to\infty,
\end{gather}
and expand the right-hand side of \eqref{eq:Jacbeta13}, denoted by $W(x)$, at $x=\ell_{n-k+1}$ for small values of~$\eps$:
\begin{gather*}
W (x_k )=W (\ell_{n-k+1}+\eps )=\sum_{k=0}^\infty \frac{\eps^k}{k!}\frac{ {\rm d}^k}{ {\rm d}x^k}W(x)=0,
\end{gather*}
where the derivatives are evaluated at $x=\ell_{n-k+1}$. Substituting the expansion of $\eps$ and those in~\eqref{eq:Jacbeta14}, we collect equal powers of $\beta$ to obtain $\eps_k$. The first few are
\begin{gather*}
\eps_1=\frac{x v_1}{\alpha+n},\\
\eps_2=-\frac{x}{2(\alpha+n)^2}\bigl(2nu_1v_1+2\alpha u_1 v_1+\alpha v_1^2-2\alpha v_2 + 2nv_1^2-2xv_1v_1^{\prime}-xv_1^2-2nv_2-v_1^2\bigr),
\end{gather*}
where $u_k$, $v_k$ are the coefficients in \eqref{eq:Jacbeta14}, and $x=\ell_{n-k+1}$. In terms of the original variables:
\begin{gather*}
\eps_1= -\frac{x}{2}(\alpha+2n+x+1),\\
\eps_2= \frac{x}{24}\bigl(5+7\alpha^2+12\alpha+24n(1+\alpha+n)+13(1+\alpha+2n)x+4x^2\bigr),\\
\eps_3 = -\frac{x}{48}\bigl(9\alpha^3+42n\alpha^2+23x\alpha^2+21\alpha^2+42x\alpha+72n^2\alpha+84xn\alpha + 15\alpha+14x^2\alpha+72n\alpha\\
\hphantom{\eps_3 =}{}+19x+72n^2+30n+48n^3 +84n^2x+3+84nx+14x^2+28nx^2+2x^3\bigr).
\end{gather*}
The coefficients $\eps_1$ and $\eps_2$ correspond with $A$ and $B$ in Theorem 2 of \cite{Dimitrov:2016:ABJ}, where a little more general result is given.\footnote{The term $2\alpha(1+2\mu)$ in $B$ of \cite[Theorem 2]{Dimitrov:2016:ABJ} seems to be not correct; it should correspond with the term $12\alpha$ in our $\eps_2$.}

We give a numerical example. We take $n=5$, $\alpha=\frac13$ and $\beta=100$, and we use the expansion in \eqref{eq:Jacbeta23} with more and more terms $\eps_k/\beta^k$. In Table~\ref{tab:tab02} we give the relative errors. As expected, the larger zeros $z_k$ (larger $k$) are computed with somewhat higher accuracy.

\begin{remark}\label{rem:rem03}
We have used $n=5$ and given the results for all zeros. For larger values of $n$ the method still works for the larger zeros, but the results become less accurate. For example, when we take $n=25$, $\alpha=\frac13$, $\beta=100$, and 5 terms in the expansion in \eqref{eq:Jacbeta23}, the zero $z_{21}$ has a~relative accuracy $0.30\times 10^{-3}$ and the largest zero $z_{25}$ has a~relative error $0.37\times 10^{-6}$. When $n=50$, $z_{45}$ has a relative error $0.10\times 10^{-1}$ and $z_{50}$ has a relative error $0.94\times 10^{-5}$.
\end{remark}

\begin{table}[t]
\caption{Relative errors in the computation of the zeros $z_k$ of $P_n^{(\alpha,\beta)}(z)$ with $n=5$, $\alpha=\frac13$, $\beta=100$ by using the asymptotic expansion \eqref{eq:Jacbeta23} with more and more terms. We compared the results with the zeros computed by Maple with Digits = 16.}\label{tab:tab02}
$$
\begin{array}{lccccc}
k & $1$ \ {\rm term} \ & $2$ \ {\rm terms} \ & $3$ \ {\rm terms} \ & $4$ \ {\rm terms} \ & $5$ \ {\rm terms} \ \\
\hline
1 &\ 0.43\times10^{-2} \ &\ 0.49\times10^{-3} \ &\ 0.54\times10^{-4}&\ 0.60\times10^{-5}&\ 0.33\times10^{-4}\tsep{2pt}\\
2 &\ 0.14\times10^{-2} \ &\ 0.13\times10^{-3} \ &\ 0.12\times10^{-4}&\ 0.11\times10^{-5}&\ 0.23\times10^{-5}\\
3 &\ 0.46\times10^{-3} \ &\ 0.36\times10^{-4} \ &\ 0.28\times10^{-5}&\ 0.22\times10^{-6}&\ 0.68\times10^{-7}\\
4 &\ 0.14\times10^{-3} \ &\ 0.90\times10^{-5} \ &\ 0.60\times10^{-6}&\ 0.40\times10^{-7}&\ 0.60\times10^{-8}\\
5 &\ 0.24\times10^{-4} \ &\ 0.14\times10^{-5} \ &\ 0.85\times10^{-7}&\ 0.50\times10^{-8}&\ 0.30\times10^{-9}\\
\hline
\end{array}
$$
\end{table}

\section{An alternative expansion}\label{sec:Jacnlargebetaalt}
The first expansion in Section~\ref{sec:Jacplargebeta} gives the expansion derived in~\cite{Dimitrov:2016:ABJ}, although obtained by a~different method. In this section we derive another expansion which performs better and has simpler coefficients. First observe that the factor $(1-x/\beta)^n$ in front of the integral \eqref{eq:Jacbeta05} is included in the expansion in~\eqref{eq:Jacbeta06}. We have included it to obtain the same expansion as in~\cite{Dimitrov:2016:ABJ}, but there is no need to do so. Next, we can obtain simpler coefficients when we take $b=n+\beta$ as a new parameter. Because the large $\beta$ parameter is replaced by a larger parameter, this may have influence on the accuracy of the results.

With this new parameter we write the representation in \eqref{eq:Jacbeta04} in the form
\begin{gather}\label{eq:Jacalt01}
P_n^{(\alpha,\beta)}\left(1-\frac{2x}{b}\right) =\frac{(1-x/b)^n}{2\pi i}\int_\calC (s+1)^{n+\alpha}\left(1-\frac{x}{b-x}s\right)^{b}\frac{ {\rm d}s}{s^{n+1}},
\end{gather}
and we expand first in powers of $s$:
\begin{gather}\label{eq:Jacalt02}
\left(1-\frac{x}{b-x}s\right)^{b}=e^{-xs}\sum_{k=0}^\infty a_k(b,x) \xi^k s^k, \qquad \xi=x/(b-x), \qquad b=\beta+n.
\end{gather}
We find for the first coefficients
\begin{gather*}
a_0(b,x) =1,\qquad a_1(b,x) =-x,\qquad a_2(b,x) =\tfrac12\big(x^2-b\big),\\
a_3(b,x) =\tfrac{1}{6}\big(3bx-2b-x^3\big),\qquad a_4(b,x) =\tfrac{1}{24}\big(3b^2-6b+8bx-6bx^2+x^4\big),\\
a_5(b,x)=\tfrac{1}{120}\big(30bx+20b^2-20bx^2-15b^2x-10bx^3-24b-x^5\big),
\end{gather*}
and obtain the finite expansion
\begin{gather}
P_n^{(\alpha,\beta)}\left(1-\frac{2x}{b}\right)=(1-x/b)^n\sum_{k=0}^n \xi^k a_k(b,x) L_{n-k}^{(\alpha+k)}(x)\nonumber\\
\hphantom{P_n^{(\alpha,\beta)}\left(1-\frac{2x}{b}\right)}{} = (1-x/b)^n\sum_{k=0}^n (-\xi)^k a_k(b,x) \frac{{\rm d}^k}{{\rm d}x^k}L_{n}^{(\alpha)}(x).\label{eq:Jacalt04}
\end{gather}
The expansion is finite because powers $s^k$ in \eqref{eq:Jacalt02} with $k>n$ absorb the pole at $s=0$ in \eqref{eq:Jacalt01}. In addition, because
\begin{gather*}
\xi^{2k}a_{2k}(b,x)=\bigO\big(b^{-k}\big),\qquad \xi^{2k+1}a_{2k+1}(b,x)=\bigO\big(b^{-{k-1}}\big),\qquad b\to \infty,
\end{gather*}
the expansion in \eqref{eq:Jacalt04} has an asymptotic character for large values of $b=\beta+n$.

To obtain an expansion in negative powers of $b$ we expand
\begin{gather*}
\left(1-\frac{x}{b-x}s\right)^{b}=e^{-xs}\sum_{k=0}^\infty \frac{d_k(x;s)} {b^k},
\end{gather*}
and we find
\begin{gather*}
d_0(x;s) =1, \qquad d_1(x;s) =-\tfrac12sx^2(s+2),\\
d_2(x;s)=\tfrac{1}{24}sx^3\big({-}24-24s-8s^2+12xs+12xs^2+3xs^3\big).
\end{gather*}
These are the alternative ones of the $c_k(n,x;s)$ in \eqref{eq:Jacbeta09} and they are much simpler because $n$ is not explicitly included.

We write as in \eqref{eq:Jacbeta11}
\begin{gather*}
d_k(x;s)=\sum_{j=0}^{2k}d_{jk}s^j,
\end{gather*}
and obtain
\begin{gather*}
P_n^{(\alpha,\beta)}\left(1-\frac{2x}{b}\right)=(1-x/b)^n\sum_{k=0}^n\frac{\Psi_k(n,\alpha,x)}{b^k},
\end{gather*}
where
\begin{gather*}
\Psi_k(n,\alpha,x)=\frac{1}{2\pi i}\int_\calC (s+1)^{n+\alpha}e^{-xs}d_k(x;s)\frac{{\rm d}s}{s^{n+1}},
\end{gather*}
and, as in \eqref{eq:Jacbeta12},
\begin{gather*}
\Psi_k(n,\alpha,x)=\sum_{j=0}^{\min(n,2k)}d_{jk}L_{n-j}^{(\alpha+j)}(x)=\sum_{j=0}^{\min(n,2k)}d_{jk}(-1)^j\frac{{\rm d}^j}{{\rm d}x^j}L_{n}^{(\alpha)}(x).
\end{gather*}

Finally, we write
\begin{gather*}
\dsp{P_n^{(\alpha,\beta)}\left(1-\frac{2x}{b}\right)}=(1-x/b)^n W(n,\alpha,\beta,x),\nonumber\\
W(n,\alpha,\beta,x) = L_{n}^{(\alpha)}(x)Y(n,\alpha,\beta,x)+L_{n-1}^{(\alpha)}(x)Z(n,\alpha,\beta,x),
\end{gather*}
with expansions
\begin{gather}\label{eq:Jacalt13}
Y(n,\alpha,\beta,x)=\sum_{k=0}^n\frac{y_k(n,\alpha,x)}{b^k},\qquad Z(n,\alpha,\beta,x)=\sum_{k=0}^n\frac{z_k(n,\alpha,x)}{b^k}.
\end{gather}
The coefficients $y_k, z_k$ can be obtained with the summations as in \eqref{eq:Jacbeta18}, and the first few are
\begin{gather*}
y_0(n,\alpha,x)=1,\qquad z_0(n,\alpha,x)=0,\\
y_1(n,\alpha,x) = \tfrac12n(2x+\alpha+1),\qquad z_1(n,\alpha,x) = -\tfrac12(n+\alpha)(\alpha+x+1).
\end{gather*}

We give a numerical example. We take $n=10$, $\alpha=\frac13$, $x=1$, and several values of $\beta$. We use the expansions in \eqref{eq:Jacalt13} summing up till $k=k_{\max}$. In Table~\ref{tab:tab03} we give the relative errors. We compared the results with computed results by Maple with Digits = 16. Comparing the results with those in Table~\ref{tab:tab01}, we see a better performance for the alternative expansion. As observed in Remark~\ref{rem:rem02}, smaller values of $x$ give better performance as well.

\begin{table}[t]
\caption{Relative errors in the computation of $P_n^{(\alpha,\beta)}(1-2x/b)$ with $n=10$, $\alpha=\frac13$, $x=1$, $b=n+\beta$, and several values of $\beta$. We use the expansion in \eqref{eq:Jacalt13} summing up till $k=k_{\max}$.}\label{tab:tab03}
$$
\begin{array}{cccccc}
k_{\max} \to & 1 & 2&3&4 &5 \\
\beta &&&&&\\
\hline
 \hphantom{00}50 & 0.10\times 10^{-1} & 0.38\times 10^{-3} & 0.13\times 10^{-4} & 0.39\times 10^{-6} & 0.12\times 10^{-7}\tsep{2pt}\\
 \hphantom{0}100 & 0.33\times 10^{-2} & 0.67\times 10^{-4} & 0.12\times 10^{-5} & 0.20\times 10^{-7} & 0.33\times 10^{-9}\\
 \hphantom{0}500 & 0.16\times 10^{-3} & 0.72\times 10^{-6} & 0.28\times 10^{-8} & 0.10\times 10^{-10} & 0.36\times 10^{-13}\\
1000 & 0.42\times 10^{-4} & 0.94\times 10^{-7} & 0.18\times 10^{-9} & 0.34\times 10^{-12} & 0.60\times 10^{-15}\\
\hline
\end{array}
$$
\end{table}

\subsection{Expansions of the zeros}\label{sec:Jacnlargebetaaltzer}
We use the method and notation of the zeros as in Section~\ref{sec:Jacnlargebetazer}. We assume for $x_k$ an expansion of the form
\begin{gather}\label{eq:Jacalt15}
x_k=\ell_{n-k+1}+\delta,\qquad \delta\sim\sum_{k=1}^\infty \frac{\delta_k}{b^k},\qquad b= \beta+n, \qquad \beta\to\infty,
\end{gather}
and the first coefficients are
\begin{gather*}
\delta_1= -\frac{x}{2}(\alpha+x+1),\\
\delta_2= \frac{x}{24}\bigl(5+7\alpha^2+12\alpha+(13+13\alpha+2n)x+4x^2\bigr),\\
\delta_3=-\frac{x}{48}\bigl(9\alpha^3+23x\alpha^2+21\alpha^2+42x\alpha+15\alpha+6xn\alpha+14x^2\alpha +2x^3+6nx+19x\\
\hphantom{\delta_3=}{} +4nx^2+14x^2+3\bigr).
\end{gather*}

We repeat the numerical example of Section~\ref{sec:Jacnlargebetazer} and the results are shown in Table~\ref{tab:tab04}. We take $n=5$, $\alpha=\frac13$ and $\beta=100$, and we use the expansions in \eqref{eq:Jacalt15} with more and more terms~$\delta_k/b^k$. In Table~\ref{tab:tab04} we give the relative errors.
\begin{table}[t]
\caption{Relative errors in the computation of the zeros $z_k$ of $P_n^{(\alpha,\beta)}(z)$ with $n=5$, $\alpha=\frac13$, $\beta=100$ by using the asymptotic expansions in \eqref{eq:Jacalt15} with more and more terms. We compared the results with the zeros computed by Maple with Digits = 16.\label{tab:tab04}}
$$
\begin{array}{lccccc}
k & $1$ \ {\rm term} \ & $2$ \ {\rm terms} \ & $3$ \ {\rm terms} \ & $4$ \ {\rm terms} \ & $5$ \ {\rm terms} \ \\
\hline
1 &\ 0.12\times10^{-2} \ &\ 0.69\times10^{-4} \ &\ 0.37\times10^{-5}&\ 0.19\times10^{-6}&\ 0.71\times10^{-3}\tsep{2pt}\\
2 &\ 0.27\times10^{-3} \ &\ 0.10\times10^{-4} \ &\ 0.41\times10^{-6}&\ 0.16\times10^{-7}&\ 0.78\times10^{-5}\\
3 &\ 0.53\times10^{-4} \ &\ 0.15\times10^{-5} \ &\ 0.41\times10^{-7}&\ 0.12\times10^{-8}&\ 0.73\times10^{-7}\\
4 &\ 0.79\times10^{-5} \ &\ 0.14\times10^{-6} \ &\ 0.27\times10^{-8}&\ 0.52\times10^{-10}&\ 0.27\times10^{-9}\\
5 &\ 0.55\times10^{-6} \ &\ 0.57\times10^{-8} \ &\ 0.61\times10^{-10}&\ 0.67\times10^{-12}&\ 0.37\times10^{-13}\\
\hline
\end{array}
$$
\end{table}

Comparing the results with those of Table~\ref{tab:tab02}, we observe a better performance of the alternative expansion. In addition, the coefficients of the expansions in \eqref{eq:Jacalt13} are more concise than those in \eqref{eq:Jacbeta14}. This holds as well for the coefficients in the expansions of the zeros.

In Remark~\ref{rem:rem03} we have shown a few results for larger values of the degree $n$. In the present case we take again $n=25$, $\alpha=\frac13$, $\beta=100$, and 5 terms in the expansion in \eqref{eq:Jacalt15}. Then, the zero~$z_{21}$ has a relative accuracy $0.82\times 10^{-7}$ and the largest zero~$z_{25}$ has a relative error $0.40\times 10^{-15}$. When we take $n=50$ we have similar results: $z_{45}$ has a~relative error $0.50\times 10^{-7}$ and $z_{50}$ has a~relative error $0.51\times 10^{-16}$. We computed the results in these examples with Maple's parameter Digits = 32, otherwise, with Digits = 16, the error in the largest zeros would have been zero.

\subsection*{Acknowledgments}
The authors thank the referees for their helpful comments. This work was supported by {\it Ministerio de Econom\'{\i}a y Competitividad}, project MTM2015-67142-P (MINECO/FEDER, UE). NMT thanks CWI, Amsterdam, for scientific support.

\pdfbookmark[1]{References}{ref}
\LastPageEnding

\end{document}